\documentclass[11pt]{amsart}

\usepackage{amsmath,amsthm,amsfonts,amssymb,amsopn}
\newtheorem{thm}{Theorem}[section]

\newtheorem{rem}[thm]{Remark}

\begin{document}

\title[Nonexistence of self-similar blowup]{Nonexistence of self-similar blowup \\for the nonlinear Dirac equations in (1+1) dimensions}

\author{Hyungjin Huh}
\address{Department of Mathematics, Chung-Ang University, Seoul 156-756, Korea}
\email{huh@cau.ac.kr}

\author{Dmitry E. Pelinovsky}
\address{Department of Mathematics and Statistics, McMaster University, Hamilton, Ontario, Canada, L8S 4K1}
\email{dmpeli@math.mcmaster.ca}

\begin{abstract}
We address a general system of nonlinear Dirac equations in (1+1) dimensions
and prove nonexistence of classical self-similar blowup solutions in the space of bounded functions.
While this argument does not exclude the possibility of finite-time blowup, it still suggests that
smooth solutions to the nonlinear Dirac equations in (1+1) dimensions
do not develop self-similar singularities in a finite time.
In the particular case of the cubic Dirac equations, we characterize (unbounded)
self-similar solutions in the closed analytical form.
\end{abstract}



\maketitle
\date{}

\section{Introduction}
\setcounter{equation}{0}

Smooth solutions of many nonlinear dispersive wave equations may blowup in a finite time depending
on the power of nonlinearity. The classical example is the nonlinear Schr\"{o}dinger equation (NLS)
in (1+1) dimensions with power nonlinearity, where smooth solutions are global only in the subcritical
case. For critical (quintic) and supercritical powers, smooth solutions of the NLS may blowup in a finite time \cite{Caz}.

Nonlinear Dirac equations are considered to be the relativistic generalization of the NLS equation, yet
they display many new dynamical properties compared to the NLS equation \cite{BC-19}. In particular,
smooth solutions to many examples of the nonlinear Dirac equations in (1+1) dimensions escape
blowup in a finite time \cite{Pel}.

The general system of massless nonlinear Dirac equations in (1+1) dimensions can be written in the form:
\begin{equation}\label{equ11}
\left\{ \begin{array}{l}
  i(\partial_t U_1 +  \partial_x U_1) = \partial_{\bar{U}_1}W(U_1,U_2,\bar{U}_1,\bar{U}_2) ,\\
  i( \partial_t U_2 -  \partial_x U_2) = \partial_{\bar{U}_2}W(U_1,U_2,\bar{U}_1,\bar{U}_2),
  \end{array} \right.
\end{equation}
where $(U_1,U_2): \Bbb{R} \times \Bbb{R} \to \Bbb{C} \times \Bbb{C}$, $\bar{U}$ is a complex conjugate of $U$,
and the nonlinear potential $W$ is assumed to satisfy the following properties:
\begin{enumerate}
\item Symmetry: $W(U_1,U_2,\bar{U}_1,\bar{U}_2) = W(U_2,U_1,\bar{U}_2,\bar{U}_1)$
\item Phase invariance: $W(e^{i \theta}U_1,e^{i \theta}U_2,e^{-i \theta} \bar{U}_1,e^{-i \theta} \bar{U}_2) =
W(U_1,U_2,\bar{U}_1,\bar{U}_2)$, $\theta \in \Bbb{R}$.
\item Homogeneous polynomial in  $(U_1,U_2,\bar{U}_1, \,\bar{U}_2)$.
\end{enumerate}
It was shown in \cite{CP} that the nonlinear potential $W$ can be characterized as a homogeneous polynomial
in variables $(|U_1|^2+|U_2|^2)$, $|U_1|^2 |U_2|^2$, and $(\bar{U}_2 U_1 + U_2 \bar{U}_1)$.
In particular, the most general quartic polynomial for $W$ is represented by
\begin{align*}
W = a_1 |U_1|^2|U_2|^2+ a_2 (\bar{U}_1U_2 + \bar{U}_2U_1)^2+a_3 (|U_1|^4+ |U_2|^4)+ a_4 (|U_1|^2+ |U_2|^2)(\bar{U}_1U_2 + \bar{U}_2 U_1),
\end{align*}
where $(a_1,a_2,a_3,a_4)$ are real constants.

When $W=  |U_1|^2|U_2|^2$,  the system \eqref{equ11} is called Thirring model \cite{Thirring}.
The Cauchy problem for the Thirring model was found to be globally well-posed
in Sobolev space $H^s(\Bbb{R})$ with $s \in \mathbb{N}$ \cite{Del} and in $L^2(\Bbb{R})$
\cite{Can, Huh, Sel}. Orbital stability of solitary wave solutions in the massive Thirring model
was proven in \cite{PS2,PS}.

When $W=(\bar{U}_1U_2 + \bar{U}_2U_1)^2$, the system \eqref{equ11} is called Gross-Neveu model \cite{GN}.
The Cauchy problem for the Gross--Neveu model was proven to be globally well-posed in $H^s(\Bbb{R})$ with $s>1/2$
\cite{Huh3,Zhang} by obtaining bounds on the $L^{\infty}(\Bbb{R})$ of the solution
and in $L^2(\Bbb{R})$ \cite{H_M,ZhangZhao} by using characteristics. Spectral stability of solitary wave
solutions in the massive Gross--Neveu model was studied numerically in the general case \cite{BerComech,Cuevas,Lakoba}
and analytically in the nonrelativistic limit \cite{BC1,BC2}.

When $W = |U_1|^4 + 4 |U_1|^2 |U_2|^2 + |U_2|^4$, the system \eqref{equ11} is called the coupled-mode model
\cite{Good}. The Cauchy problem for the coupled-mode system was found to be well-posed in Sobolev space $H^s(\Bbb{R})$
with $s \in \mathbb{N}$ \cite{Good,Pel} and in $L^2(\Bbb{R})$ \cite{Huh2}. Existence and spectral stability
of solitary wave solutions have been analyzed in this model in many details (see \cite{CP} and references therein).

Finally, when $W=(|U_1|^2+  |U_2|^2) (\bar{U}_1U_2 + \bar{U}_2U_1)$, the nonlinear Dirac equation
with pseudoscalar potential \cite{Stu} occurs in the context of photonic crystals with the nonlinear refractive index \cite{AP}.
As far as we know, it has been an open problem for many years to address global existence   or finite time blowup of solutions to the Cauchy problem for this system \cite{Pel}. This problem is the subject of the present work.

The self-similar blowup has played an important role in the  formation of singularities of 
partial differential equations (see \cite{EF} and references therein). In particular, 
self-similar blowup solutions  have been investigated for the relativistic wave equation \cite{B,K,MZ} 
and for the Navier--Stokes equations \cite{Tsai}.

The main goal of this study is to prove nonexistence of classical self-similar blowup
solutions in the space of bounded functions to the nonlinear Dirac equation (\ref{equ11})
with the nonlinear potential in the form:
\begin{equation}
\label{W}
W = (|U_1|^2+ |U_2|^2)^k (\bar{U}_1U_2 + \bar{U}_2 U_1)^{\ell},
\end{equation}
where $k, \,\ell$ are nonnegative integers with $p := k + \ell - 1 \in \mathbb{N}$.
Besides the space and time translation invariance, the system of nonlinear Dirac equations
\eqref{equ11} with (\ref{W}) has the following scaling invariance property:
if $[U_1(x,t), \,U_2(x,t)]$ is a solution, then
$$
\left[\lambda^{\frac{1}{2p}} U_1 (\lambda x, \lambda t), \,\lambda^{\frac{1}{2p}} U_2 (\lambda x, \lambda t)\right], \quad
\lambda > 0
$$
is also a solution of the same system.
Thanks to the scaling invariance property and the separation of variables,
the class of self-similar solutions is defined in the form:
\begin{equation}
\label{self-similar}
U_1(x,t) = \frac{1}{(1-t)^{\frac{1}{2p}}} U\left( \frac{x}{1-t}\right), \quad
U_2(x,t) = \frac{1}{(1-t)^{\frac{1}{2p}}} V\left( \frac{x}{1-t}\right),
\end{equation}
where $U$ and $V$ are functions of $y := x/(1-t)$. A singularity of the self-similar solutions (\ref{self-similar}) 
is placed at the point $(x,t) = (0,1)$ thanks to the space and time translation
symmetries. Thanks to the unit speed of propagation,
the variable $y$ can be restricted to the interval $[-1,1]$.

Existence of classical solutions $(U,V) : [-1,1] \mapsto \mathbb{C} \times \mathbb{C}$
implies self-similar blowup of smooth solutions to the Cauchy problem to the nonlinear Dirac
equations (\ref{equ11}) in a finite time.
While nonexistence of classical solutions does not exclude
the possibility of finite-time blowup, it still suggeists that
smooth solutions to the nonlinear Dirac equations
do not develop self-similar singularities in a finite time.

The following theorem presents the main result of this work.

\begin{thm}
\label{theorem-1}
For every $p \in \mathbb{N}$, there exist no self-similar solutions
in the form \eqref{self-similar} with $(U,V) \in C^1(-1,1) \cap L^{\infty}([-1,1])$.
\end{thm}

Theorem \ref{theorem-1} is proven in Section 2 by using the polar decomposition,
dynamical system methods, and a continuation argument. The proof is simpler
in the case of odd $\ell$ and more technically involved in the case of even $\ell$.

The main reason for nonexistence of classical self-similar solutions is 
their breakdown either before they reach the end points $y = \pm 1$ 
of the interval $[-1,1]$ or at the end points $y = \pm 1$. In the general case,
we are not able to obtain the precise rate of how $U$ and $V$ diverge before or at $y = \pm 1$.
However, in the case $k = \ell = 1$, which corresponds
to the physically relevant model (\ref{equ11}) with $W=(|U_1|^2+  |U_2|^2) (\bar{U}_1U_2 + \bar{U}_2U_1)$ 
derived in \cite{AP},
we are able to integrate the system of differential equations for $U$ and $V$
provided the initial condition satisfies $|U(0)| = |V(0)|$.

The following theorem represents the result.

\begin{thm}
\label{theorem-2}
For $k = \ell = 1$, there exists a unique classical self-similar solution in the form \eqref{self-similar}
with $|U(0)| = |V(0)|$ that extends to $y \to 1$ and satisfies the following asymptotic behavior
\begin{equation}
\label{asympt-scaling}
U(y) \sim (1-y)^{\frac{1}{4}}, \quad V(y) \sim (1-y)^{-\frac{1}{4}} \quad \mbox{\rm as} \quad y \to 1.
\end{equation}
However, this solution does not extend to $y \to -1$ in the sense that
there exists $y_0 \in (-1,0)$ such that $\lim_{y \to y_0} U(y)$ and $\lim_{y \to y_0} V(y)$ diverge.
All other classical solutions with $|U(0)| = |V(0)|$ extend neither to $y \to 1$ nor to $y \to -1$.
\end{thm}

Theorem \ref{theorem-2} is proven in Section 3, where the system of differential equations for $U$ and $V$
with the initial condition $|U(0)| = |V(0)|$ is integrated in a closed form. It follows from the same method
used in the proof of Theorem \ref{theorem-2} (see Remark \ref{rem-blow} below), that all other solutions
with $|U(0)| \neq |V(0)|$ do not extend simultaneously to $y \to 1$ and $y \to -1$. Therefore,
the classical self-similar solutions do not exist in Theorem \ref{theorem-1} because they blow up
before reaching $y = \pm 1$ at least for $k = \ell = 1$.

\section{Proof of Theorem \ref{theorem-1}}
\setcounter{equation}{0}

Substituting \eqref{self-similar} into \eqref{equ11} with $W$ in (\ref{W}) yields
the following system of differential equations for $U$ and $V$:
\begin{align}\label{equ_UV}
\begin{cases}
\, i \big[(y+1) U' +\frac{1}{2p} U  \big] = F V + G U, \\
\, i\big[(y-1) V' +\frac{1}{2p} V  \big] = F U + G V,
\end{cases}
\end{align}
where the prime denotes derivative in $y$, $p = k + \ell - 1$, and
$$
F = \ell (|U|^2+|V|^2 )^k ( U \bar{V} + \bar{U} V)^{\ell-1}, \qquad
G = k (|U|^2+|V|^2 )^{k-1} ( U \bar{V} + \bar{U} V)^{\ell}.
$$
We are studying existence of classical solutions to the system (\ref{equ_UV})
on the interval $[-1,1]$, hence $(U,V) \in C^1(-1,1)$. We also require
the classical solution to remain bounded as $y \to \pm 1$, hence
$(U,V) \in L^{\infty}([-1,1])$.

By inspecting the integrating factors for the left-hand side
of the system (\ref{equ_UV}), we introduce the new variables:
\begin{equation}
u(y) := (1+y)^{1/2p} U(y), \qquad v(y) := (1-y)^{1/2p} V(y).
\label{new-variables-u-v}
\end{equation}
New variables allow us to rewrite the system \eqref{equ_UV} in the equivalent form:
\begin{equation}
\label{equ_uv}
\left\{ \begin{array}{l}
i (1+y)^{1-\frac{1}{2p}} u' = F v (1-y)^{-\frac{1}{2p}} + G u (1+y)^{-\frac{1}{2p}}, \\
-i (1-y)^{1-\frac{1}{2p}} v' = F u (1+y)^{-\frac{1}{2p}} + G v (1-y)^{-\frac{1}{2p}},
\end{array} \right.
\end{equation}
where
\begin{eqnarray*}
F = \ell \left( \frac{|u|^2}{(1+y)^{\frac{1}{p}}} + \frac{|v|^2}{(1-y)^{\frac{1}{p}}} \right)^k
\frac{(u \bar{v} + \bar{u} v)^{\ell-1}}{(1-y^2)^{\frac{\ell - 1}{2p}}}
\end{eqnarray*}
and
\begin{eqnarray*}
G = k \left(\frac{|u|^2}{(1+y)^{\frac{1}{p}}} + \frac{|v|^2}{(1-y)^{\frac{1}{p}}} \right)^{k-1}
\frac{(u \bar{v} + \bar{u} v)^{\ell}}{(1-y^2)^{\frac{\ell}{2p}}}.
\end{eqnarray*}
If $(U,V) \in C^1(-1,1) \cap L^{\infty}([-1,1])$, then $(u,v) \in C^1(-1,1) \cap L^{\infty}([-1,1])$
and moreover, $u(-1) = 0$ and $v(1) = 0$.

Let us use the polar decomposition for complex-valued amplitudes:
\begin{equation}
\label{polar-form}
u = |u| e^{i \alpha}, \quad v = |v| e^{i \beta},
\end{equation}
where all functions depend on $y$. Because $(u,v) \in C^1(-1,1)$, then
$$
\frac{d|u|}{dy}, \quad \frac{d|v|}{dy}, \quad |u| \frac{d \alpha}{d y}, \quad |v| \frac{d \beta}{dy}
$$
are all bounded and piecewise continuous on the interval $(-1,1)$. Therefore, substituting the polar decomposition (\ref{polar-form})
into (\ref{equ_uv}) and separating the real and imaginary parts yield the following system of differential equations
for amplitudes and phases:
\begin{align}\label{equa}
\left\{
\begin{aligned}
\frac{d|u|}{dy} & = \frac{F |v| \sin(\beta - \alpha)}{(1+y)^{1-\frac{1}{2p}}  (1-y)^{\frac{1}{2p}}},\\
\frac{d|v|}{dy} & = \frac{F |u| \sin(\beta - \alpha)}{(1+y)^{\frac{1}{2p}}  (1-y)^{1-\frac{1}{2p} }},
\end{aligned}\right.
\end{align}
and
\begin{align}\label{equ-phase}
\left\{
\begin{aligned}
-|u|\frac{d \alpha}{dy} & = \frac{F |v| \cos(\beta - \alpha)}{(1+y)^{1-\frac{1}{2p}}  (1-y)^{\frac{1}{2p}}} +
\frac{G |u|}{(1+y)},\\
|v|\frac{d \beta}{dy} & = \frac{F |u| \cos(\beta - \alpha)}{(1+y)^{\frac{1}{2p}} (1-y)^{1-\frac{1}{2p}}} +
\frac{G |v|}{(1-y)},
\end{aligned}\right.
\end{align}
where
\begin{eqnarray*}
F = \ell \left( \frac{|u|^2}{(1+y)^{\frac{1}{p}}} + \frac{|v|^2}{(1-y)^{\frac{1}{p}}} \right)^k
\frac{[2|u||v| \cos(\beta - \alpha)]^{\ell-1}}{(1-y^2)^{\frac{\ell - 1}{2p}}}
\end{eqnarray*}
and
\begin{eqnarray*}
G = k \left(\frac{|u|^2}{(1+y)^{\frac{1}{p}}} + \frac{|v|^2}{(1-y)^{\frac{1}{p}}} \right)^{k-1}
\frac{[2|u||v| \cos(\beta - \alpha)]^{\ell}}{(1-y^2)^{\frac{\ell}{2p}}}.
\end{eqnarray*}
The vector field of the system (\ref{equa}) and (\ref{equ-phase}) is piecewise continuous on $(-1,1)$.
We shall now proceed differently depending whether $\ell$ is zero, odd, or even.

\subsection{The case of $\ell = 0$}

In this case, $F = 0$ and the system (\ref{equa}) implies that $|u(y)|$ and $|v(y)|$ are constant in $y$.
Therefore, it is impossible to satisfy $u(-1) = 0$ and $v(1) = 0$ except for the trivial (zero) solution.

 \subsection{The case of odd $\ell$}

Since $\ell$ is odd, we have $F \geq 0$ and $G \cos(\beta - \alpha) \geq 0$.
Combining the two equations in the system (\ref{equ-phase}) yields
\begin{eqnarray}
\nonumber
|u||v| \frac{d}{dy} \sin (\beta-\alpha)
 & = & F \cos^2(\beta-\alpha)
 \left( \frac{|v|^2}{(1+y)^{1-\frac{1}{2p}} (1-y)^{\frac{1}{2p}}}
 + \frac{|u|^2}{(1+y)^{\frac{1}{2p}} (1-y)^{1-\frac{1}{2p}}} \right) \\
 & & \quad + \frac{2 G |u||v| \cos(\beta - \alpha)}{1-y^2} \geq 0.
 \label{monotonicity}
\end{eqnarray}
From here, we obtain a contradiction against the existence of
solutions $(u,v) \in C^1(-1,1) \cap L^{\infty}([-1,1])$
satisfying $u(-1) = 0$ and $v(1) = 0$.

Indeed, if $u(-1) = 0$, then $\frac{d}{dy} |u| \geq 0$ at least near $y = -1$.
The first equation of the system (\ref{equa}) with odd $\ell$
implies $\sin(\beta - \alpha) \geq 0$ at least near $y = -1$.
Thanks to monotonicity (\ref{monotonicity}),
we have $\sin(\beta - \alpha) \geq 0$ for every $y \in (-1,1)$.
The second equation of the system (\ref{equa}) with odd $\ell$ implies
then that $\frac{d}{dy} |v| \geq 0$ for every $y \in (-1,1)$.
Hence $|v(y)| \geq |v(-1)|$ for every $y \in (-1,1)$ and
it is impossible to satisfy $v(1)=0$ except for the trivial (zero) solution.

 \subsection{The case of even $\ell$}

Since $\ell$ is even, we have $G \geq 0$ and $F \cos(\beta - \alpha) \geq 0$.
Combining the two equations in the system (\ref{equ-phase}) yields
\begin{eqnarray}
\nonumber
|u||v| \frac{d}{dy} (\beta-\alpha)
 & = & F \cos(\beta-\alpha)
 \left( \frac{|v|^2}{(1+y)^{1-\frac{1}{2p}} (1-y)^{\frac{1}{2p}}}
 + \frac{|u|^2}{(1+y)^{\frac{1}{2p}} (1-y)^{1-\frac{1}{2p}}} \right) \\
 & & \quad + \frac{2 G |u||v|}{1-y^2} \geq 0.
 \label{monotonicity-even}
\end{eqnarray}
If $\ell \geq 2$, then $F = G = 0$ if $\cos(\beta - \alpha) = 0$, hence
$\beta - \alpha = \pm \frac{\pi}{2}$ are invariant lines, which cannot be
crossed for finite $y \in (-1,1)$. From here, we obtain a contradiction against the existence of
solutions $(u,v) \in C^1(-1,1) \cap L^{\infty}([-1,1])$
satisfying $u(-1) = 0$ and $v(1) = 0$.

Indeed, if $u(-1) = 0$, then $\frac{d}{dy} |u| \geq 0$ at least near $y = -1$.
The first equation of the system (\ref{equa}) with even $\ell$
implies $\sin(\beta - \alpha) \cos(\beta - \alpha) \geq 0$ at least near $y = -1$.
Thanks to monotonicity (\ref{monotonicity-even}) and invariance of $\beta - \alpha = \pm \frac{\pi}{2}$,
we have
\begin{align*}
\mbox{either} \quad 0 \leq \beta-\alpha \leq  \frac{\pi}{2} \quad
\mbox{or} \quad \pi \leq \beta-\alpha \leq  \frac{3\pi}{2},
\end{align*}
for every $y \in (-1,1)$, which means that
$\sin(\beta - \alpha) \cos(\beta - \alpha) \geq 0$ for every $y \in (-1,1)$.
The second equation of the system (\ref{equa}) with even $\ell$ implies then
that $\frac{d}{dy} |v| \geq 0$ for every $y \in (-1,1)$.
Hence $|v(y)| \geq |v(-1)|$ for every $y \in (-1,1)$ and
it is impossible to satisfy $v(1)=0$ except for the trivial (zero) solution.

\section{Proof of Theorem \ref{theorem-2}}
\setcounter{equation}{0}

Here we investigate the case of $k = \ell = 1$ in the system \eqref{equa} and \eqref{equ-phase}.
The system is rewritten explicitly as follows:
\begin{align}\label{u-v}
\left\{
\begin{aligned}
\frac{d|u|}{dy} & = \frac{|v| \sin(\beta - \alpha)}{\sqrt{1-y^2}} \left[ \frac{|u|^2}{1+y} + \frac{|v|^2}{1-y} \right],\\
\frac{d|v|}{dy} & = \frac{|u| \sin(\beta - \alpha)}{\sqrt{1-y^2}} \left[ \frac{|u|^2}{1+y} + \frac{|v|^2}{1-y} \right],
\end{aligned}
\right.
\end{align}
and
\begin{align}
\label{alpha-beta}
\left\{
\begin{aligned}
-|u|\frac{d \alpha}{dy} & = \frac{|v| \cos(\beta - \alpha)}{\sqrt{1-y^2}}\left[ \frac{3 |u|^2}{1+y} + \frac{|v|^2}{1-y} \right],\\
|v|\frac{d \beta}{dy} & = \frac{|u| \cos(\beta - \alpha)}{\sqrt{1-y^2}} \left[ \frac{|u|^2}{1+y} + \frac{3 |v|^2}{1-y} \right].
\end{aligned}
\right.
\end{align}
We are looking for classical solutions $(u,v) \in C^1(-1,1)$ satisfying
the constraint $|u(0)| = |v(0)|$ on the initial condition.
The system (\ref{u-v}) yields the first-order invariant
\begin{align} \label{rela2}
|u(y)|^2 =|v(y)|^2 + C,
\end{align}
where $C$ is constant. It follows from the constraint $|u(0)| = |v(0)|$ that $C = 0$, hence
$|u(y)| = |v(y)|$ for every $y \in [-1,1]$.
With this reduction, the system (\ref{u-v}) and (\ref{alpha-beta}) reduces to a simpler form:
\begin{equation}
\label{u-alpha}
\left\{
\begin{aligned}
\frac{d|v|}{dy} & = \frac{2 |v|^3 \sin(\beta - \alpha)}{\sqrt{(1-y^2)^3}}, \\
\frac{d (\beta-\alpha)}{dy} & = \frac{8 |v|^2 \cos(\beta - \alpha)}{\sqrt{(1-y^2)^3}}.
\end{aligned}
\right.
\end{equation}

Let us introduce the independent variable $\tau : [-1,1] \mapsto \mathbb{R}$ by
\begin{equation}
\label{tau}
\tau(y) := \int_0^y \frac{dy}{\sqrt{(1-y^2)^3}}.
\end{equation}
Then, $\tau(y) \to \pm \infty$ as $y \to \pm 1$. Let us also rewrite the system (\ref{u-alpha})
in dependent variables
\begin{equation}
\label{xi-eta}
\xi := |v|, \quad \eta := \sin(\beta - \alpha).
\end{equation}
Then, the system (\ref{u-alpha}) can be written as the autonomous planar dynamical system:
\begin{equation}
\label{u-dynam}
\left\{
\begin{aligned}
\dot{\xi} & = 2 \xi^3 \eta, \\
\dot{\eta} & = 8 \xi^2 (1-\eta^2),
\end{aligned}
\right.
\end{equation}
where the dot denotes derivative with respect to $\tau$. The line segment 
$\Sigma_0 := \{ \xi = 0, \; \eta \in [-1,1]\}$ consists of the degenerate critical
points, whereas $\Sigma_{\pm} := \{ \xi \in \mathbb{R}, \; \eta = \pm 1\}$ 
are invariant lines with the one-dimensional flow in $\xi$ given by $\dot{\xi} = \pm 2 \xi^3$.

The system (\ref{u-dynam}) is integrable
with the first invariant $E(\xi,\eta) := \xi^8 (1-\eta^2)$, where the values of $E$ are constant
and $E \geq 0$ since $\eta \in [-1,1]$. The value $E = 0$ is not isolated since
$\Sigma_0$ intersects $\Sigma_{\pm}$. The flow on $\Sigma_{\pm}$ is given by $\dot{\xi} = \pm 2 \xi^3$.
For both signs, $\xi(\tau)$ does not exist for every $\tau \in \mathbb{R}$ since it blow-up in a finite $\tau$
either before $\tau \to -\infty$ or before $\tau \to +\infty$.
For the minus sign, the solution satisfies $\lim_{\tau \to +\infty} \xi(\tau) = 0$ (that is,
$|u(y)| \to 0$ as $y \to 1$) and moreover
$$
\xi(\tau) \sim \tau^{-1/2} \quad \Rightarrow \quad |u(y)| \sim (1-y)^{1/4}.
$$
This provides the asymptotic scaling (\ref{asympt-scaling}) in variables $U$ and $V$
thanks to the transformation (\ref{new-variables-u-v}).

For every $E > 0$, the level curve $E(\xi,\eta) = E > 0$ is unbounded in $\xi$ and does not intersect 
$\Sigma_0$ or $\Sigma_{\pm}$. It follows from the second equation in the system (\ref{u-dynam})
that $\dot{\eta} > 0$, hence the map $\tau \to \eta$ is strictly monotonic along the flow
with $\eta \in (-1,1)$. It follows from the first equation of the system (\ref{u-dynam})
by the comparison principle with $\eta \geq \eta_0 > 0$ that the map $\tau \to \xi$ blows up
before $\tau \to +\infty$ in positive flow in $\tau$. Similarly, by the comparison principle with $\eta \leq -\eta_0 < 0$
it follows that the map $\tau \to \xi$ blows up before $\tau \to -\infty$ in negative flow in $\tau$.
Therefore, no other solutions bounded near $y = 1$ exist.

\begin{rem}
\label{rem-blow}
For general initial conditions with $|u(0)| \neq |v(0)|$, we have $C \neq 0$
in the local invariant \eqref{rela2}. For solutions bounded near $y = 1$, we have $v(1) = 0$
and $C = |u(1)|^2 > 0$. For solutions bounded near $y = -1$, we have $u(-1) = 0$
and $C = -|v(-1)|^2 < 0$. In the former case $C > 0$, the system
of differential equations takes the form
\begin{equation}
\label{u-alpha-C}
\left\{
\begin{aligned}
\frac{d|v|}{dy} & = \frac{\sqrt{C + |v|^2} \sin(\beta - \alpha)}{\sqrt{(1-y^2)^3}} \left[ 2 |v|^2 + C(1-y) \right], \\
\frac{d (\beta-\alpha)}{dy} & = \frac{\cos(\beta - \alpha)}{\sqrt{(1-y^2)^3} |v| \sqrt{C + |v|^2}} \left[ 8 |v|^2 + 2C (4-y) |v|^2 + C^2 (1-y) \right],
\end{aligned}
\right.
\end{equation}
so that the same definitions for $\tau$, $\xi$ and $\eta$ as in \eqref{tau} and \eqref{xi-eta} can be employed.
The same monotonicity argument for the map $\tau \to \eta$ and the same comparison principle
for the map $\tau \to \xi$ can be employed to show that the solutions bounded near $y = 1$ blow up at a finite $y_0 \in (-1,1)$ before
the other end $y = -1$. However, if $C \neq 0$, solutions extending to $y \to 1$ do not satisfy
the same asymptotic behaviour as in \eqref{asympt-scaling}.
\end{rem}

\vspace{0.25cm}

{\bf Acknowledgment:} H. Huh was supported by LG Yonam Foundation of Korea.
D. Pelinovsky was supported by the NSERC Discovery grant.

 \end{document}